\documentclass[12pt,leqno]{article}
\catcode`\@=11 \long\def\@makecaption#1#2{
 \vskip 10pt
 \setbox\@tempboxa\hbox{\bf #1: \sf #2}
 \ifdim \wd\@tempboxa >\hsize \bf #1: \sf #2\par \else \hbox
to\hsize{\hfil\box\@tempboxa\hfil}
 \fi}
\catcode`\@=12
\usepackage{amsmath}
\usepackage{amsfonts}
\usepackage{amssymb}
\usepackage{epsfig}
\newcommand{\al}{\alpha}
\newcommand{\be}{\beta}
\newcommand{\R}{\mathbb R}

\newcommand{\SSS}{\mathbb S}

\newcommand{\bbx}{\hfill\rule{.46em}{.46em}}
\newcommand{\inter}[1]{{\rm int}\,#1}

\newcommand{\conv}{{\rm conv}}

\newcommand{\rd}[1]{\R^{#1}}

\newtheorem{theo}{Theorem}
\newtheorem{lem}[theo]{Lemma}

\newtheorem{cor}[theo]{Corollary}
\newtheorem{pr}[theo]{Proposition}

\begin{document}
\title{Shadow systems and volumes of polar convex bodies
\footnotetext{2000 Mathematics Subject Classification 52A40.}
}
\author{Mathieu Meyer\footnote{Both authors have been supported in part
by the Research Network Program in Mathematics of the High Council for
Scientific and Technological Cooperation between France and Israel.}  \\
{\small Equipe d'Analyse et de Math\'ematiques Appliqu\'ees}\\
{\small Universit\'e de Marne-la-Vall\'ee}\\
{\small Cit\'e Descartes}\\
{\small 5, boulevard Descartes}\\
{\small Champs-sur-Marne}\\
{\small 77454 Marne-la-Vall\'ee}\\
{\small Cedex 2, France}\\
{\small \tt meyer@math.univ-mlv.fr}\\
\and
\setcounter{footnote}{0}
Shlomo Reisner\footnotemark\hspace{1mm} \footnote{Shlomo Reisner has been supported
in part by NATO Collaborative Linkage Grant PST.CLG 979701.}\\
{\small Department of Mathematics}\\
{\small University of Haifa}\\
{\small Haifa 31905, Israel}\\
{\small \tt reisner@math.haifa.ac.il}}
\date{}
\maketitle
\begin{abstract}
We prove that the reciprocal of the volume of the polar bodies,
about the Santal\'o point, of a {\em shadow system\/} of convex bodies
$K_t$, is a convex function of $t$. Thus extending to the non-symmetric
case a result of Campi and Gronchi. The case that the reciprocal of the
volume is an affine function of $t$ is also investigated and is
characterized under certain conditions.

We apply these results to prove exact reverse Santal\'o inequality for
polytopes in $\rd{d}$ that have at most $d+3$ vertices.
\end{abstract}

\section{Introduction and notations}

A {\em convex body\/} in $\rd{d}$ is a compact convex set with non-empty
interior. If $K$ is a convex body in $\rd{d}$ and $z$ is an interior point
of $K$, the {\em polar body\/} of $K$ with respect to $z$, $K^{*z}$, is
defined by
\[K^{*z}=\{y\in \rd{d}\,;\, \forall x\in K,\;\;\langle y,x-z\rangle\leq 1\}\]
($\langle\cdot,\cdot\rangle$ is the canonical scalar product in $\rd{d}$).
Note that in some other places the polar body of $K$ with respect to $z$ is
defined as a translation by $z$ of the above, namely:
$\{y\in \rd{d}\,;\, \forall x\in K,\;\;\langle y-z,x-z\rangle\leq 1\}$.
We denote by $|A|$ the $k$-dimensional volume of a
measurable set $A\subset \rd{d}$, where $k$ is the dimension of the minimal
flat containing $A$ ({\em volume} means $k$-dimensional Lebesgue measure in
this flat).
A well known result of Santal\'o \cite{Sant} states that in every convex body
$K\subset \rd{d}$ there exists a unique point $S(K)$ - the {\em Santal\'o
point\/} of $K$, such that
\[|K^{*S(K)}|=\min_{z\in \inter(K)}|K^{*z}|\,,\]
 We shall denote $K^{*S(K)}$ by $K^*$.

A {\em shadow system\/} along a direction $\theta \in \SSS^{d-1}$ is a
family of convex
sets $K_t\subset \rd{d}$ which are defined by
\[K_t=\conv\{x+\alpha(x)t\theta\,;\,x\in M\subset \rd{d}\}\,,\]
where $M$ is a bounded set, $\alpha$ is a bounded function on $M$ and $t$
belongs to an interval in $\R$ ($\conv(A)$ is the closed convex hull of a set
$A\subset \rd{d}$). We say that the shadow system $K_t$ is
{\em non-degenerate\/} if $K_t$ has non-empty interior for all $t$ in the
interval.

Shadow systems were introduced by Rogers and
Shephard \cite{R-S} in order to treat extremal problems for convex bodies.
They proved that $t\mapsto |K_t|$ is a convex function of $t$. The concept was
further investigated by Shephard \cite{Sh} who, among other results, extended
the convexity result to mixed volumes.

Campi and Gronchi proved in a recent paper \cite{C-G}: {\em If $K_t$ is a
shadow system of origin symmetric convex bodies in $\rd{d}$, then
$|K^*_t|^{-1}$ is a convex function of\ \ \!$t$.} They applied this result to
prove reverse forms of the $L^p$-Blaschke-Santal\'o inequality of Lutwak and
Zhang \cite{L-Z} in dimension 2 (as well as to provide a new proof of the
result of \cite{L-Z}).

In Section~2 of this paper we prove a result (Theorem~\ref{th-A}) analogous
to Campi and Gronchi's result, avoiding the symmetry assumption. The proof
in this, more general, setting is more delicate and requires developing
other methods. We also
investigate the case that $t\mapsto |K^*_t|^{-1}$ is affine and prove that
if $t\mapsto |K_t|$ is affine then $t\mapsto |K^*_t|^{-1}$ is affine only
if all the bodies $K_t$ are affine images of each other.

In Section~3 we apply the results of Section~2 to prove ``exact reverse
Santal\'o inequality'' for polytopes in $\rd{d}$ that have at most $d+3$
vertices.

A well known conjecture, called sometimes ``Mahler's conjecture'', states
that, for every convex body $K$ in $\rd{d}$,
\begin{equation}
\label{eq-{-A}}
\Pi_d(K)=|K||K^*|\geq \Pi_d(\Delta)=\frac{(d+1)^{d+1}}{(d!)^2}\,,
\end{equation}
where $\Delta$ is a $d$-dimensional simplex ($\Pi_d(K)$ is called the
{\em volume-product\/} of $K$). It is also conjectured that
equality in (\ref{eq-{-A}}) is attained {\em only} if $K$ is a simplex.
The inequality (\ref{eq-{-A}}) for $d=2$ was proved by Mahler \cite{Mah}
with the case of equality proved by Meyer \cite{Mey-2}. Other cases, like
e.g.\ bodies of revolution, were treated in \cite{Mey-Rei}. Several special
cases in the centrally symmetric situation can be found in
\cite{SR,Rei-1,GMR,Mey-1,Rei-2}. The (non-exact) reverse Santal\'o
inequality of Bourgain and Milman \cite{B-M} is
\[\Pi_d(K)\geq c^d\Pi_d(B)\]
where $c$ is a positive constant and $B$ is the Euclidean ball (or any
ellipsoid). This should be compared with the Blaschke-Santal\'o inequality
\[\Pi_d(K)\leq \Pi_d(B)\]
with equality only for ellipsoids (\cite{Sant}, \cite{Pet}, see \cite{M-P}
for a simple proof of both the inequality and the case of equality).

We prove in this paper that if $K$ is a convex polytope in $\rd{d}$ with
at most $d+3$ vertices, then
$\Pi_d(K)$ is never less then $\Pi_d(\Delta)$, where $\Delta$ is
a $d$-dimensional simplex. Equality holds only if $K$ itself is a simplex.

In the last section we present, as another application of the tools
developed in Section~2, new proofs of two well known theorems. One is
Blaschke-Santal\'o inequality. The proof presented here looks smooth, in
particular
the characterization of the maximal bodies as ellipsoids is a natural
part of the proof of the inequality itself and is simpler than in other
known proofs.

The second is reverse Santal\'o inequality in dimension 2 (non-symmetric
case) \cite{Mah}, together with the characterization of the minimal
bodies as triangles \cite{Mey-2}. Again, the new method provides a unified
and simple proof of the inequality together with the case of equality.

The notations we use are standard notations of the theory of convex bodies,
as may be found e.g.\ in R. Schneider's book \cite{Sch}. We refer the
reader to this book for background material as well. In particular, a
convex {\em polytope\/} is a convex body which is the convex hull of
finitely many points (vertices) in $\rd{d}$. A {\em pyramid\/} is a
convex body which is the convex hull in $\rd{d}$ of a point (apex) with
a $(d-1)$-dimensional convex body (basis). A {\em double pyramid\/} is
the convex hull of a $(d-1)$-dimensional convex body $F$ in $\rd{d}$, and
two points $x_1,x_2$, such that $x_1$ and $x_2$ are on opposite sides of
the hyperplane containing $F$, and the line segment $]x_1,x_2[$ intersects
$F$.

\section{The convexity of $t\mapsto |K^*_t|^{-1}$}

\begin{theo}
\label{th-A}
Let $K_t$, $t\in [a,b]$ be a non-degenerate shadow system in $\rd{d}$.
Then $|K^*_t|^{-1}$ is a convex function of $t\in [a,b]$.
\end{theo}

The inequality stated in the following lemma is a particular case of
a result due to K. Ball \cite{Ball} (see also \cite{Bus}). We present
here a proof taken from
\cite{M-P} in order to specify the conditions for equality.

\begin{lem}
\label{lem-A}
Let $f,g,h:\R_+\rightarrow \R_+$ be functions which are
compactly supported and continuous on their supports.
Assume further that for all $y,z>0$
\begin{equation}
\label{eq-A}
f\left({2zy\over z+y}\right) \ge g(y)^{z\over z+y}h(z)^{y\over z+y}
\end{equation}
Then
$${1\over \int_0^{\infty} f(t)\,dt}\le {1\over 2}
\left({1\over\int_0^{\infty} g(t)\,dt}+
{1\over \int_0^{\infty} h(t)\,dt}\right)$$
with equality if and only if denoting $B=\int g$ and $C=\int h$,
we have for all $x\geq 0$
$$g(Bx)= h(Cx)=f\left({2BCx\over B+C}\right)\,.$$
\end{lem}
\vskip 1mm\noindent
{\bf Proof. } For $u\in [0,1]$, we define the functions $y(u)$ and
$z(u)$ by
$$\int_{0}^{y(u)} g(t) dt =Bu\hbox{ and } \int_{0}^{z(u)} g(t) dt
=Cu.$$
One gets $$y'(u)={B\over g(y(u))}\hbox{ and } z'(u)={C\over h(z(u))}\,.$$
Now, setting $$x(u)={2y(u)z(u)\over y(u)+z(u)}\, \hbox{ which implies }\,
x'=2{z^{2}y'+y^2 z'\over (y+z)^2}$$ we get
$$\int f(x) dx\ge  2\int_{0}^{1}g(y)^{z\over z+y}h(z)^{z\over z+y}
\left({z^2 B\over g(y)}+ {y^{2} C\over h(z)}\right) {1\over (y+z)^2} du\,.
$$
Since
$${1\over z+y}\left({z^2 B\over g(y)}+ {y^{2} C\over h(z)}\right)=
{z\over z+y}{z B\over g(y)}+{y\over z+y}{y C\over h(z)}\ge
\left({z B\over g(y)}\right)^{z\over z+y}\left({y C\over h(z)}\right)^{y\over
z+y}\, $$
with equality if and only if ${z B\over g(y)}={y C\over h(z)}$,  it follows that
$$\int f(x) dx\ge 2\int_{0}^1 (zB)^{z\over z+y} (yC)^{y\over z+y}
{1\over z+y} du.$$
Setting $\lambda= {z\over z+y}$, we have
$$2(zB)^{z\over z+y} (yC)^{y\over z+y}{1\over
z+y}= {2BC\over \left({C\over \lambda}\right)^{\lambda}
\left({B\over 1-\lambda}\right)^{1-\lambda}}\ge {2BC\over B+C},$$
with equality if and only if ${C\over \lambda}={B\over 1-\lambda}$,
that is $Bz=Cy$.
Thus the result follows, with equality if and only if for every $u$,
$Bz(u)=Cy(u)$, $g(y(u))=h(z(u))=f\left({2y(u)z(u)\over z(u)+y(u)}\right)$.
This means, setting $x(u)={y(u)\over B}={z(u)\over C}$, that
$$g(Bx)=h(Cx)=f\left( {2BCx\over B+C} \right)$$ for every $x\ge 0$.\bbx

Without loss of generality,
we may and shall assume throughout this section,  that the shift vector
$\theta$ from the definition of a shadow body
is the $d$-th coordinate unit vector of $\rd{d}$. That is, representing
$\rd{d}$ as $\rd{d-1}\times \R$, we have for all $(X,x)\in M$ ($M\subset
\rd{d}$ being a bounded set) a velocity $v(X,x)$ such that $K_t$ is the
closed convex hull of $\{(X,x+tv(X,x))\,;\,(X,x)\in M\}$. We denote by $P$
the orthogonal projection of $\rd{d}$ onto $\rd{d-1}$. For a convex body
$K\subset \rd{d}$ and for $y\in \R$ we denote
\[K(\cdot,y)=\{Y\in \rd{d-1}\,;\,(Y,y)\in K\}\,,\]
and similarly, for $Y\in \rd{d-1}$,
\[K(Y,\cdot)=\{y\in \R\,;\,(Y,y)\in K\}\,.\]

\begin{lem}
\label{lem-B}
Let $C\in\rd{d-1}$ be an interior point of $\conv(P(M))$. For
$a\leq s<t\leq b$ let $a_s$ and $a_t$ be interior points of $K_s$ and
$K_t$ respectively. Let $a_{\frac{s+t}{2}}=\frac{a_s+a_t}{2}$ and assume
that $a_{\frac{s+t}{2}}\in \inter K_{\frac{s+t}{2}}$.
Let $G_u=(C,a_u)$, $u=s,t\hbox{ or }\frac{s+t}{2}$.
Define $g(y)=|K_s^{*G_s}(\cdot,y)|$, $h(z)=|K_t^{*G_t}(\cdot,z)|$,
$f(x)=|K_{\frac{s+t}{2}}^{*G_{\frac{s+t}{2}}}(\cdot,x)|$. Then the functions
$g$, $h$ and $f$ satisfy the assumptions of Lemma~\ref{lem-A} for all $y,z>0$.
\end{lem}
\vskip 1mm\noindent
{\bf Proof. } For an interior point $G=(C,\alpha)$ of $K_t$
one may describe the polar body $K_t^{*G}$ of $K_t$ with
respect to $G$ in the following way:
$$K_t^{*G}=\{ (Y,y); \langle X-C,Y\rangle+(x+v(X,x)t-\alpha)y\le 1\hbox{
for every $(X,x)\in M$} \}.$$
\vskip 1mm\noindent
For $y,z> 0$ let $Y\in K_s^{*G_s}(\cdot,y)$ and $Z\in
K_t^{*G_t}(\cdot,z)$.
Then for $(X,x)\in M$ we have
$$\langle X-C,zY+yZ\rangle
 + 2 \left(x+v(X,x){s+t\over 2}- {a_s+a_t\over 2}\right)zy\leq y+z\,.$$
\vskip 1mm\noindent
Therefore we get for every $z,y>0$,
$${zY+yZ\over z+y}\in K_{s+t\over 2}^{*G_{s+t\over 2}}\left(\cdot,{2zy\over
z+y}\right), $$
that is
$${z\over z+y} K_s^{*G_s}(\cdot,y) +{y\over z+y}  K_t^{*G_t}(\cdot,z)
\subset K_{s+t\over 2}^{*G_{s+t\over 2}}\left(\cdot,{2zy\over
z+y}\right).$$
Inequality~(\ref{eq-A}) now follows  from the Brunn-Minkowski inequality.
\bbx
\vskip 3mm
As an immediate result of Lemma~\ref{lem-A} and Lemma~\ref{lem-B} we get
\begin{lem}
\label{lem-C}
With the notations of Lemma~\ref{lem-B},
let $B_+(u,a_u)=\int_0^\infty |K_u^{*G_u}(\cdot,x)|\,dx$
and $B_-(u,a_u)=\int_{-\infty}^0 |K_u^{*G_u}(\cdot,x)|\,dx$. Then
\begin{equation}
\label{eq-AA}
\frac{1}{B_+(\frac{s+t}{2},a_{\frac{s+t}{2}})}\leq \frac{1}{2}\left(
\frac{1}{B_+(s,a_s)}+\frac{1}{B_+(t,a_t)}\right)\,.
\end{equation}
The same inequality holds for $B_-$ instead of $B_+$.
\end{lem}

One can verify now that for every $X\in \rd{d-1}$
\begin{equation}
\label{eq-AAA}
K_{\frac{s+t}{2}}(X,\cdot)\subset \frac{1}{2}(K_s(X,\cdot)+K_t(X,\cdot))
\end{equation}
(which is how the result of \cite{R-S} on the convexity of $|K_t|$ is proved).
If for every $u$ such that $K_u(C,\cdot)\neq\emptyset$ we denote
$K_u(C,\cdot)=[\alpha_u,\beta_u]$ then (\ref{eq-AAA}) means
\begin{equation}
\label{eq-B}
\frac{\alpha_s+\alpha_t}{2}\leq \alpha_{\frac{s+t}{2}}<\beta_{\frac{s+t}{2}}
\leq \frac{\beta_s+\beta_t}{2}\,.
\end{equation}
\begin{lem}
\label{lem-D}
With the above notations, and those of the preceding lemmas,
for every $a\in ]\al_{\frac{s+t}{2}},
\be_{\frac{s+t}{2}}[$ there exist $a_s\in ]\al_s,\be_s[,\,\,
a_t\in ]\al_t,\be_t[$, such that $a=a_{\frac{s+t}{2}}=
\frac{a_s+a_t}{2}$ and
\[\frac{B_+(s,a_s)}{B_-(s,a_s)}=\frac{B_+(t,a_t)}{B_-(t,a_t)}\,.\]
\end{lem}
\vskip 2mm\noindent
{\bf Proof. }For $v\in \R$ and $u=s$ or $t$, let $G_u=(C,v)$,
$B_+(u,v)=\int_0^{\infty}|K_u^{*G_u}(\cdot,x)|\,dx$ and
$B_-(u,v)=\int_{-\infty}^0|K_u^{*G_u}(\cdot,x)|\,dx$.
The functions $v\mapsto B_{+}(s,v)$ and $v\mapsto
B_{+}(t,v)$ are continuous on the intervals
$]\alpha_{s},\beta_{s}[$ and
$]\alpha_{t},\beta_{t}[$ respectively. They are bounded from below
by positive numbers
and they tend to $+\infty$ on the right-hand side of these intervals.
Similarly, the functions $v\mapsto B_{-}(s,v)$ and $v\mapsto
B_{-}(t,v)$ are continuous, are bounded from below by positive numbers
on the same intervals, and tend to infinity
on their left-hand sides.

Define $\rho: ]\max(\alpha_{s},
2a-\beta_{t}),\min(\beta_{s},2a-\alpha_{t}[\rightarrow \R$ by
$$\rho(v)={B_{+}(s,v)\over B_{-}(s,v)}-{B_{+}(t,2a-v)\over B_{-}(t,2a-v)}$$
(note that the assumption on $a$, together with (\ref{eq-B}), imply that the open
interval of definition of $\rho$ is not empty).

Now
\vskip 1mm\noindent
- If $ 2a-\beta_{t}<\alpha_{s}$ then ${B_{+}(s,v)\over B_{-}(s,v)}\rightarrow 0$
and ${B_{+}(t,2a-v)\over B_{-}(t,2a-v)}$ is bounded from below by a positive
number as $v\rightarrow \alpha_{s}$.
\vskip 1mm\noindent
- If  $ 2a-\beta_{t}\geq\alpha_{s}$ then ${B_{+}(s,v)\over B_{-}(s,v)}$ is bounded
from above
and ${B_{+}(t,2a-v)\over B_{-}(t,2u-w)}\rightarrow \infty$ as $v\rightarrow
2a-\beta_{t}$.
\vskip 1mm\noindent
So $\rho$ is negative on the left of its
interval of definition. In a similar way we see that it is  positive on
its right. It follows from continuity that
$\rho$ vanishes at some point $v=a_s$ in the open interval. Defining
$a_t=2a-a_s$ we get the result. \bbx
\vskip 2mm\noindent
{\bf Proof of Theorem~\ref{th-A}. }We want to prove that
\begin{equation}
\label{eq-C}
\frac{1}{|K_{\frac{s+t}{2}}^*|}\leq \frac{1}{2}\left(\frac{1}{|K_s^*|}+
\frac{1}{|K_t^*|}\right)\,.
\end{equation}

Let $(C,a)$ be the Santal\'o point of $K_{\frac{s+t}{2}}$. By Lemma~\ref{lem-D}
there are points $a_u\in \inter K_u(C,\cdot)$, $u=s$ or $t$, such that
$a=a_{\frac{s+t}{2}}=\frac{a_s+a_t}{2}$ and
\[\frac{B_+(s,a_s)}{B_-(s,a_s)}=\frac{B_+(t,a_t)}{B_-(t,a_t)}=\lambda\]
If we denote again $G_u=(C,a_u)$, $u=s$,  $t$ or $\frac{s+t}{2}$, we have
$|K_u^{*G_u}|=B_+(u,a_u)+B_-(u,a_u)$. Thus, for $u=s$ or $t$ we have
\[B_-(u,a_u)=\frac{1}{1+\lambda}|K_u^{*G_u}|\,.\]

We conclude now from Lemma~\ref{lem-C} that
\[|K_{\frac{s+t}{2}}^*|=|K_{\frac{s+t}{2}}^{*G_{\frac{s+t}{2}}}|=
B_+(\frac{s+t}{2},a_{\frac{s+t}{2}})+B_-(\frac{s+t}{2},a_{\frac{s+t}{2}})
\geq\]
\[\frac{2B_+(s,a_s)B_+(t,a_t)}{B_+(s,a_s)+B_+(t,a_t)}+
\frac{2B_-(s,a_s)B_-(t,a_t)}{B_-(s,a_s)+B_-(t,a_t)}=\]
\[(1+\lambda)\frac{2B_-(s,a_s)B_-(t,a_t)}{B_-(s,a_s)+B_-(t,a_t)}=
\frac{2|K_s^{*G_s}||K_t^{*G_t}|}{|K_s^{*G_s}|+|K_t^{*G_t}|}\,.\]

Thus
\begin{equation}
\label{eq-D}
\frac{1}{|K_{\frac{s+t}{2}}^*|}\leq \frac{1}{2}\left(
\frac{1}{|K_s^{*G_s}|}+\frac{1}{|K_t^{*G_t}|}\right)\leq
\frac{1}{2}\left(\frac{1}{|K_s^{*}|}+\frac{1}{|K_t^{*}|}\right)\,,
\end{equation}
where the very last inequality is due to the minimality of $|K_u^*|$.
\bbx

\subsection{The case of equality}
Let $a\leq s<t\leq b$ and assume that equality holds in the inequality
(\ref{eq-C}). Then, from (\ref{eq-D}) and uniqueness of the Santal\'o
point, it follows that, in the notations of the proof of Theorem~\ref{th-A},
$G_s$ and $G_t$ are, respectively, the Santal\'o points of $K_s$ and $K_t$.
In other words, the Santal\'o points of $K_s$, $K_t$ and $K_{\frac{s+t}{2}}$
are, respectively, $(C,a_s)$, $(C,a_t)$, $(C,{\frac{a_s+a_t}{2}})$.

Moreover, from the convexity of $|K_u^*|^{-1}$, proven in Theorem~\ref{th-A},
it follows that if equality holds in (\ref{eq-C}) then for every
$u=(1-\al)s+\al t\in ]s,t[$, ($0<\al<1$), we have
\[|K_u^*|^{-1}=(1-\al)|K_s^*|^{-1}+\al |K_t^*|^{-1}\,.\]
This in its turn implies, again by the argument of the proof of
Theorem~\ref{th-A}, that the Santal\'o point of $K_u$ is
$(C,(1-\al)a_s+\al a_t)$. Furthermore, the way of proof of Theorem~\ref{th-A}
shows that for all such $u$
\[\frac{B_+(u,a_u)}{B_-(u,a_u)}\equiv \lambda\,,\]
the same $\lambda$ for all $u\in ]s,t[$. Here $a_u=(1-\al)a_s+\al a_t$.

We summarize these facts in the following lemma and add to them one more fact,
which is a consequence of Lemmas \ref{lem-A} and \ref{lem-B}.

\begin{lem}
\label{lem-E}
Let $a\leq s<t\leq b$. If equality holds in the inequality
$|K_{\frac{s+t}{2}}^*|^{-1}\leq \frac{1}{2}(|K_s^*|^{-1}+
|K_t^*|^{-1})$ then:
\begin{itemize}
\item[1)] $u\mapsto |K_u^*|^{-1}$ is affine on the interval $[s,t]$.
\item[2)] The Santal\'o points of $K_u$ have a fixed orthogonal projection $C$
on $\rd{d-1}$ and they behave affinely on $[s,t]$. That is,
$S(K_u)=(C,(1-\al)a_s+\al a_t)$, for $u=(1-\al)s+\al t$ and $0<\alpha<1$.
\item[3)] The ratio $\lambda=\frac{B_+(u,a_u)}{B_-(u,a_u)}$, $a_u=(1-\al)a_s+\al a_t$
is fixed for $u=(1-\al)s+\al t\in [s,t]$.
\item[4)] There exists a function $\phi:\,J\rightarrow \rd{d-1}$ such that
for all $u\in [s,t]$ we have for all $x\in J=\{x\in \R\,;\,
K_{\frac{s+t}{2}}^*(\cdot,|K_{\frac{s+t}{2}}^*|x)\neq\emptyset\}$
\begin{equation}
\label{eq-E}
K_u^*(\cdot,|K_u^*|x)=K_{\frac{s+t}{2}}^*(\cdot,|K_{\frac{s+t}{2}}^*|x)+
\left(u-\frac{s+t}{2}\right)|K_u^*|\phi(x)\,.
\end{equation}
\end{itemize}
\end{lem}

\vskip 1mm\noindent
{\bf Proof. }We have only to prove 4). Let us consider first $x>0$. Going back
to the proofs of Lemmas \ref{lem-B} and \ref{lem-C}, we see that the properties
1) and 2) above show that if equality holds then $s$ and $t$ in
Lemmas \ref{lem-B} and \ref{lem-C} may be replaced by any $u,v\in [s,t]$, $u<v$
(provided that $G_u=(C,a_u)$ is the Santal\'o point of $K_u$). Using the
equality case in Brunn-Minkowski inequality at the conclusion of the proof of
Lemma~\ref{lem-B}, and the characterization of equality in Lemma~\ref{lem-A},
we conclude that for every $u\in [s,t]$ and $x\geq 0$,
$K_u^*(\cdot,B_+(u,a_u)x)$ is a translate of
$K_{\frac{s+t}{2}}^*(\cdot,B_+({\frac{s+t}{2}},a_{\frac{s+t}{2}})x)$.
Say for some function $\mu:\R\times
[s,t]\rightarrow \rd{d-1}$, we have
$$K_u^{*}(\cdot,B_+(u,a_u)x) =
K_{s+t\over 2}^{*}\left(\cdot,B_+\left(\frac{s+t}{2},a_{\frac{s+t}{2}}\right)
x\right)+\mu(x,u)\,.$$
>From the equality (which is the equality case of the inclusion at the end
of the proof of Lemma~\ref{lem-B})
\[\frac{B_+(v,a_v)}{B_+(u,a_u)+B_+(v,a_v))}K_u^{*}(\cdot,B_+(u,a_u)x) +
\frac{B_+(u,a_u)}{B_+(u,a_u)+B_+(v,a_v))}K_v^{*}(\cdot,B_+(v,a_v)x)=\]
\[K_{\frac{u+v}{2}}^{*}\left(\cdot,B_+\left({\frac{u+v}{2}},
a_{\frac{u+v}{2}}\right)x\right)\,,\]
we conclude that
\[\frac{B_+(v,a_v)}{B_+(u,a_u)+B_+(v,a_v))}\mu(x,u)+
\frac{B_+(u,a_u)}{B_+(u,a_u)+B_+(v,a_v))}\mu(x,v)=
\mu\left(x,\frac{u+v}{2}\right)\]
which means, since by 1) and 3) above
$u\mapsto B_+(u,a_u)^{-1}$ is affine on $[s,t]$,
that for every $x\geq 0$,
$u\mapsto \mu(x,u)B_+(u,a_u)^{-1}$ is an affine function of
$u$. Since, by the
definition, $\mu(x,\frac{s+t}{2})=0$, it follows that for $x\geq 0$
\[K_u^*(\cdot,B_+(u,a_u)x)=K_{\frac{s+t}{2}}^*\left(\cdot,
B_+\left(\frac{s+t}{2},a_{\frac{s+t}{2}}\right)x\right)+
\left(u-\frac{s+t}{2}\right)B_+(u,a_u)\phi_+(x)\,.\]
for some function $\phi_+:\,J_+\rightarrow \rd{d-1}$, where $J_+$
is an appropriate interval. In a similar way we get
for $x\leq 0$
\[K_u^*(\cdot,B_-(u,a_u)x)=K_{\frac{s+t}{2}}^*\left(\cdot,
B_-\left(\frac{s+t}{2},a_{\frac{s+t}{2}}\right)x\right)+
\left(u-\frac{s+t}{2}\right)B_-(u,a_u)\phi_-(x)\,.\]

We substitute, in view of 3) above,
$B_+(v,a_v)=\frac{\lambda}{1+\lambda}|K_v^*|$ and
$B_-(v,a_v)=\frac{1}{1+\lambda}|K_v^*|$ in the above equalities.
Making the required changes of variables and defining $\phi$
accordingly, separately for $x\geq 0$ and $x\leq 0$, we get the
equality (\ref{eq-E}).\vspace{.25em} \bbx

Using Lemma~\ref{lem-E} we get the following result.

\begin{pr}
\label{pr-A}
Let $K_t$, $t\in [a,b]$, be a non-degenerate shadow system in
$\rd{d}$. Then the following are equivalent:
\begin{itemize}
\item[1)] $t\mapsto |K_t|$ and
$t\mapsto |K^*_t|^{-1}$ are both affine functions of
$t\in [a,b]$.
\item[2)] There exit real numbers $v$ and $u$, and a vector
$V\in \rd{d-1}$, such that, for all
$t\in [a,b]$, $K_t$ is the image of $K_{\frac{a+b}{2}}$
under the affine transformation $A_t:\,\rd{d}\rightarrow \rd{d}$
which, when $\rd{d}$ is represented as $\rd{d-1}\times\R$, is
given by
\begin{equation}
\label{eq-EE}
A_t(X,x)=\left(X,x+\left(t-\frac{a+b}{2}\right)(vx+\langle
V,X\rangle+u)\right)\,.
\end{equation}
\end{itemize}
\end{pr}

\vskip 1mm\noindent
{\bf Problem. } Can the assumption on $|K_t|$,  in the
``$1)\Rightarrow 2)$''
direction of Proposition~\ref{pr-A}, be replaced by
a weaker one while keeping the conclusion that the bodies $K_t$
must be affine images of each other true?

\vskip 1mm\noindent
{\bf Proof. }We first prove that 1) implies 2). Assume that
both $t\mapsto |K_t|$ and $t\mapsto |K_t^*|^{-1}$ are affine.

One may assume, for convenience, that $[a,b]=
[-c,c]$ is a symmetric interval. We first make a couple of
observations.

For $X\in \rd{d-1}$ let $K_t(X,\cdot)=[a_t(X),b_t(X)]$.
By the hypothesis of the proposition we have
\[\left|K_{\frac{s+t}{2}}\right|=\frac{1}{2}|K_s|+
\frac{1}{2}|K_t|\]
that is
\[\int_{\rd{d-1}}\left(b_{\frac{s+t}{2}}(X)-
a_{\frac{s+t}{2}}(X)\right)\,dX=
\int_{\rd{d-1}}\left[\frac{1}{2}(b_s(X)-a_s(X))+
\frac{1}{2}(b_t(X)-a_t(X))\right]\,dX\,.\]
But by (\ref{eq-AAA}) we have
\[[a_{\frac{s+t}{2}}(X),b_{\frac{s+t}{2}}(X)]\subset
\frac{1}{2}[a_s(X),b_s(X)]+\frac{1}{2}[a_t(X),b_t(X)]\,.\]
We conclude that the support of $a_t(X)$ and $b_t(X)$
is $P(K_0)$ and, for every $X\in P(K_0)$,
$t\mapsto a_t(X)$ and $t\mapsto b_t(X)$ are affine
functions. It follows that there exist functions
$p_1,p_2:P(K_0)\rightarrow \R$, such that for all
$t\in [-c,c]$
\[a_t(X)=a_0(X)+tp_1(X)\,,\quad b_t(X)=b_0(X)+tp_2(X)\,\]
(this last fact is actually equivalent to $t\mapsto |K_t|$
being affine).

By Lemma~\ref{lem-E}, 2), all the Santal\'o points of $K_t$,
$t\in [-c,c]$, have the same projection on $\rd{d-1}$. Moreover,
by an appropriate translation, we may assume that the Santal\'o
point of $K_0$ is $(0,0)$. Then, again by Lemma~\ref{lem-E}, the
Santal\'o points of $K_{-c}$ and $K_{c}$ are $(0,-\alpha)$ and
$(0,\alpha)$, respectively. We define a new shadow system
$\tilde K_t$ by changing the speed function from $v(X,x)$ to
$v(X,x)-\frac{\alpha}{c}$. Then, for every $t$, $\tilde K_t$ is
a translation of $K_t$ by $-\frac{\alpha}{c}t$. Thus the Santal\'o
points of $\tilde K_{-c}$ and $\tilde K_c$ are both $(0,0)$. By
Lemma~\ref{lem-E}, $(0,0)$ the Santal\'o point of
$\tilde K_t$ for all $t\in [-c,c]$. The bottom line of all
this is that we may assume, without loss of generality, that
the Santal\'o points of all the bodies $K_t$ are $0\in \rd{d}$.
We shall assume also that $|K_0^*|=1$.

Let us denote $|K_t^*|=c_t$, then $c_0=1$ and, by the assumptions
of the proposition, $c_t^{-1}=\gamma t+1$ for some constant $\gamma$.

Lemma~\ref{lem-E}, 4) implies that
\begin{equation}
\label{eq-F}
K_t^*(\cdot,c_tx)=K_0^*(\cdot,x)+tc_t\phi(x)
\end{equation}
for some function $\phi:\,J\rightarrow \rd{d-1}$ ($J$ the
appropriate interval).

\vskip 1mm
We fix now $Z\in \rd{d-1}$ with Euclidean norm $1$, and consider,
for $t\in [-c,c]$ the $2$-dimensional body
$L_{t}:=\{(x,y)\ ;\ (xZ,y)\in K_{t}\}$.
Then we have
\[L_t=\{(x,y)\,;\,xZ\in P(K_0), y\in
[a_0(xZ)+tp_1(xZ),b_0(xZ)+tp_2(xZ)]\}\,.\]

Since the polar of a plane section through the center of polarity
of a convex body is the orthogonal projection of the polar body on the
same plane,
setting $q(y')=\langle Z,\phi(y') \rangle$ for $y'$ in the interval of
definition of $\phi$, we have
$$L_{t}^{*0}=\{\big(x'+tc_tq(y'),c_ty'\big); (x',y')\in L^*\}\,.$$

It follows from the next lemma that $q(x)=kx$  for every $x$ where $q$
is defined. Thus, $x\rightarrow
\langle Z,\phi(x)\rangle$ is linear on the segment where it is defined.
It follows that $x\rightarrow \phi(x)$ is linear, that is
$\phi(x)=xV$ for some fixed vector $V\in \R^{d-1}$.
Thus
\[K^{*}_{t}= H_t(K^{*}_{0})\,,\mbox{\ where\ }H_t(X,x)= (X+tc_txV,c_tx)\,.\]
that is, $K^{*}_{t}$ is an affine image of $K_0^{*}$. It follows that $K_{t}$
is an affine image of $K_0$:
\[K_t=H^{*-1}_t(K_0)\,,\mbox{\ where\ }H^{*-1}_t(X,x)=(X,c_t^{-1}x-t\langle
X,V\rangle)\,.\]

\begin{lem}
\label{lem-F}
Suppose that $I$ and $J$ are intervals
and $a,b:I\rightarrow \R$, $\alpha,\beta:J\rightarrow \R$ are four
functions, with $a\le b$, $\alpha\leq \beta$, such that
$$L=\{(x,y); x\in I, y\in [a(x),b(x)]\}$$ is a convex body
containing $0$ in its interior and
$$L^*= \{(x',y'); y'\in J, x'\in
[\alpha(y'),\beta(y')]\}\,,$$
where $L^*$ is the polar, about $0$, of $L$.
Suppose that for some functions
$p_1,p_2:I\rightarrow \R$ and $t$ in an interval $[-c,c]$,  the set
$$L_{t}=\{x,y); x\in I, y\in [a(x)+tp_1(x),b(x)+tp_2(x)]\},$$
is convex.  We suppose, moreover,
that for some function $q: J\rightarrow\R$, and for $t\in [-c,c]$,
$$L_{t}^{*}=\{(x',c_ty'); y'\in J, x'\in
[\alpha(y')+tc_tq(y'),\beta(y')+tc_tq(y')]\}\,,$$
where $c_t=(\gamma t+1)^{-1}$ for some constant $\gamma$
(polarity is taken here again about $0$).

Then for some constant
$k\in \R$, one has for every $x\in I$ and $y'\in J$,
\begin{equation}
\label{eq-FF}
q(y')=ky'\,,\mbox{\ }p_1(x)=\gamma a(x)-kx
\mbox{\ and\ }p_2(x)=\gamma b(x)-kx\,.
\end{equation}
\end{lem}

\vskip 1mm \noindent
{\bf Proof of Lemma~\ref{lem-F}.}
\vskip 1mm\noindent Let $X=(x,y)\in L$
and $X'=(x',y')\in L^{*}$ satisfy $\langle X,X'\rangle =xx'+yy'=1$,
then $X=(x,b(x))$ or $X=(x,a(x))$ and
$X'=(\alpha(y'), y')$ or $X'=(\beta(y'),y')$.
In the case that $x\ge 0$ and $y'\ge 0$, we have $X=(x,b(x))$ and
$X'=(\beta(y'),y')$.
It follows from the hypotheses that in this case
\begin{equation}
\label{eq-FF1}
x \beta(y')+y'b(x)=1\mbox{ and}
\end{equation}
$$x(\beta(y')+tc_t q(y'))+c_t y'(b(x)+tp_2(x))\le 1$$
for all $t$.
It follows that
$$tc_t (x q(y') +y'p_2(x))+ y'b(x)(c_t-1)\le 0$$
or, dividing by $c_t$ (which is positive),
$$t (x q(y') +y'p_2(x))+ y'b(x)(1-{1\over c_t})\le 0\,.$$

Since ${c_t}^{-1}= \gamma t+1$ one gets
$$t (x q(y') +y'p_2(x)-\gamma y'b(x))\le 0$$
for all $t\in [-c,c]$. It follows that
\begin{equation}
\label{eq-FF2}
x q(y') +y'p_2(x)-\gamma y'b(x)=0\,.
\end{equation}

(\ref{eq-FF1}) and (\ref{eq-FF2}) imply that, in fact,
\begin{equation}
\label{eq-FF3}
x(\beta(y')+tc_t q(y'))+c_t y'(b(x)+tp_2(x))= 1\,.
\end{equation}

Now,  (\ref{eq-FF1}) and (\ref{eq-FF3}) show that vectors normal
to the convex bodies $L_t$ at their boundary points
$(x,b(x)+tp_2(x))$ are $(\beta(y')+tc_tq(y'),c_ty')$. It follows
that, at points where both the derivatives $b'(x)$ and $p'_2(x)$
exist (which, by convexity, are all, but at most countably many,
points $x\in I\cap \R_+$), we have
\begin{equation}
\label{eq-FF4} b'(x)=-\frac{\beta(y')}{y'} \mbox{ and }
b'(x)+tp'_2(x)= \frac{-(\beta(y')+tc_tq(y'))}{c_ty'}\,.
\end{equation}

Combining together the two equalities in (\ref{eq-FF4}), we get
\[b'(x)\big(1-\frac{1}{c_t}\big)+tp'_2(x)=-t\frac{q(y')}{y'}\]
or
\begin{equation}
\label{eq-FF5}
-\gamma b'(x)+p'_2(x)=-\frac{q(y')}{y'}\,.
\end{equation}

Now, (\ref{eq-FF2}) can be written as
\[-\frac{q(y')}{y'}=\frac{-\gamma b(x)+p_2(x)}{x}\,,\]
thus, the continuous (and Lipschitz on closed intervals contained in $I$)
function $f(x)=-\gamma b(x)+p_2(x)$, has derivatives at all points of $I$,
except, possibly, at a countable set, and satisfies the differential
equation
\[\frac{f'(x)}{f(x)}=\frac{1}{x}\]
at all points $0<x\in I$ of differentiability. By a standard argument it
follows that
\[ -\gamma b(x)+p_2(x)=-k_1x, \mbox{ that is } p_2(x)=\gamma b(x)-k_1x\]
for some constant $k_1$. By (\ref{eq-FF5}) we get
\[q(y')=k_1y'\,.\]
(There are various kinds of arguments that can be applied to overcome the
possible lack of differentiability of $f$ at {\em all} points of
$I\cap \R_+$. One such argument can make use,  e.g., of the
{\em Corollaire} at the end of \S 2 in Ch.\ 1 of \cite{Bourb}.)

Considering in a similar way the three other
cases for the signs of $x$ and $y'$, we get
\[x>0,\,y'>0\;\Rightarrow\;p_2(x)=\gamma b(x)-k_1x,\,q(y')=k_1y'\,,\]
\[x<0,\,y'>0\;\Rightarrow\;p_1(x)=\gamma a(x)-k_2x,\,q(y')=k_2y'\,,\]
\[x>0,\,y'<0\;\Rightarrow\;p_2(x)=\gamma b(x)-k_3x,\,q(y')=k_3y'\,,\]
\[x<0,\,y'<0\;\Rightarrow\;p_1(x)=\gamma a(x)-k_4x,\,q(y')=k_4y'\,.\]
Finally, these four conditions show that all the constants $k_i$ are
equal to one constant $k$ and that (\ref{eq-FF}) is satisfied. One can
then check the correctness of the representations, given above, of the
linear maps $H_t$ and $H^{*-1}_t$, by considering the relation between
the matrices $\alpha_t$ and $\alpha_t^{*-1}$ in the sequel.

This completes the proofs of Lemma~\ref{lem-F}. \bbx

The ``$2)\Rightarrow 1)$'' direction of Proposition~\ref{pr-A} is simpler.
Substitute $s=t-\frac{a+b}{2}$.
 The matrix representing the linear part of $A_t$ is
\begin{equation*}
\alpha_t=\left(
\begin{matrix}
1&\cdots&0&0 \\
\vdots&\ddots&\vdots&\vdots \\
0&\cdots&1&0\\
sv_1&\cdots&sv_{d-1}&(vs+1)
\end{matrix}
\right)\,.
\end{equation*}
The transformation that maps $K_{\frac{a+b}{2}}^*$ onto
$K_t^*$ is then represented by
\begin{equation*}
\alpha_t^{*-1}=\left(
\begin{matrix}
1&\cdots&0&-sv_1(vs+1)^{-1}\\
\vdots&\ddots&\vdots&\vdots  \\
0&\cdots&1&-sv_{d-1}(vs+1)^{-1}\\
0&\cdots&0&(vs+1)^{-1}
\end{matrix}
\right)\,.
\end{equation*}

That is,
\[|K_t|=|K_{\frac{a+b}{2}}|\det \alpha_t=|K_{\frac{a+b}{2}}|(vs+1)\,.\]
and
\[|K_t^*|=|K_{\frac{a+b}{2}}^*|\det \alpha_t^{*-1}=|K_{\frac{a+b}{2}}^*|
(vs+1)^{-1}\,.\]
Thus $t\mapsto |K_t|$ and  $t\mapsto |K_t^*|^{-1}$ are both affine and,
moreover, $|K_t||K_t^*|$ is constant for $t\in [a,b]$.
\bbx
\medskip

\noindent{\bf Remark.\ }The details of the proof of the direction
``$1)\Rightarrow 2)$'' of Proposition~\ref{pr-A} imply the following
consequence, which we use in the sequel: {\em Let
\[K_t=\{(X,x)\in \rd{d}\times \R\,;\,X\in \conv(P(M))\,,\,
a_t(X)\leq x\leq b_t(X)\}\]
and assume that $K_t$ satisfies the condition 1) of Proposition~\ref{pr-A}.
Then the transformation $A_t$, given in (\ref{eq-EE}) satisfies
\[(X,a_t(X))=A_t(X,a_{\frac{a+b}{2}}(X))\mbox{\ and\ }
(X,b_t(X))=A_t(X,b_{\frac{a+b}{2}}(X))\]
for all $X\in \conv(P(M))$.\/}
\medskip

The ``$2)\Rightarrow 1)$'' direction in the next corollary is always true,
while the ``$1)\Rightarrow 2)$'' direction is a consequence of
Proposition~\ref{pr-A}.

\begin{cor}
\label{cor-A}
Let $K_t$ $t\in [a,b]$, be a non-degenerate shadow system in
$\rd{d}$, and assume that the function
$t\mapsto |K_t|$ is affine. Then the following are equivalent:
\begin{itemize}
\item[1)] $\Pi_d(K_t)$ is constant for $t\in [a,b]$.
\item[2)] For every $s,t\in [a,b]$, $K_t$ and $K_s$ are affine images of
one another.
\end{itemize}
\end{cor}

\section{Reverse Santal\'o inequality for polytopes with few vertices}

As an application of the results of Section~2, we prove here that the
conjectured ``exact reverse Santal\'o inequality'' $\Pi_d(K)\geq
\frac{(d+1)^{d+1}}{(d!)^2}$ for convex bodies $K$ in $\rd{d}$, is valid
if $K$ is a polytope with few vertices.

\begin{theo}
\label{th-B}
Let $K$ be a convex body in $\rd{d}$, which is a polytope with at most
$d+3$ vertices. Let $\Pi_d(K)=|K||K^*|$. Then
\[\Pi_d(K)\geq \frac{(d+1)^{d+1}}{(d!)^2}\,,\]
with equality if and only if $K$ is a $d$-dimensional simplex.
\end{theo}

For the proof of Theorem~\ref{th-B} we need the following lemma.

\begin{lem}
\label{lem-G}
Let $K=\conv(\{x_0\}\cup F)$ be a pyramid in $\rd{d}$, where $F$ is a
$(d-1)$-dimensional convex body and $x_0$ is not in the hyperplane
containing $F$. Then
\[\Pi_d(K)=\frac{(d+1)^{d+1}}{d^{d+2}}\Pi_{d-1}(F)\]
and the Santal\'o point $S(K)$ of $K$ lies on the line segment $]z_0,x_0[$,
where $z_0$ is the Santal\'o point of $F$ (considered as a
$(d-1)$-dimensional convex body) and
\[\frac{|x_0-z_0|}{|S(K)-z_0|}=d+1\,.\]
\end{lem}

\vskip 1mm \noindent
{\bf Proof. }Writing $\rd{d}$ as $\rd{d-1}\times \R$, we may assume
(using an affine transformation) that $S(K)=(0,0)$, $x_0=(0,1-\alpha)$,
and that $F$ lies in the hyperplane $\{(Y,-\alpha)\,;\,Y\in \rd{d-1}\}$,
for some $0<\alpha<1$. Thus, $(0,-\alpha)$ is in the relative interior
of $F$. It is easy to check that $K^*=\conv(\{(0,-\alpha^{-1})\},\tilde F)$
where $\tilde F$ is a $(d-1)$-dimensional convex body contained in the
hyperplane $\{(Z,(1-\alpha)^{-1})\,;\,Z\in \rd{d-1}\}$, and if $F$ is
identified with $K(\cdot,-\alpha)\subset \rd{d-1}$ and $\tilde F$ with
$K(\cdot,(1-\alpha)^{-1})\subset \rd{d-1}$ then
$\tilde F=(1-\alpha)^{-1}F^{*0}$. Thus
\[|K||K^*|=|F||F^{*0}|\cdot\frac{1}{d^2\alpha}
\left(\frac{1}{1-\alpha}\right)^d\,.\]

A simple optimization shows that the minimum over $\alpha$ of the right
hand side is obtained for $\alpha=\frac{1}{d+1}$. Clearly the right hand
side is then minimal if and only if $0$ is the Santal\'o point of
$K(\cdot,-\alpha)$ (identified with $F$).\bbx

\vskip 1mm \noindent
{\bf Proof of Theorem~\ref{th-B}. }For $k\geq d+1$ let ${\cal P}_{d,k}$
be the set of convex bodies in $\rd{d}$, which are polytopes having at
most $k$ vertices. ${\cal P}_{d,k}$ is closed in the Hausdorff metric.
Also, since $\Pi_d$ is affinely invariant, we may, using F. John theorem,
restrict our attention to a compact subset of ${\cal P}_{d,k}$ to realize
that $\Pi_d$ attains a minimum in ${\cal P}_{d,k}$.
We shall call a
polytope at which a minimum is attained in ${\cal P}_{d,d+3}$, {\em minimal}.

For $d=2$ the theorem
is a very special case of the results of Mahler \cite{Mah}. We check
separately the cases of $d+2$ and $d+3$ vertices.

\vskip 1mm \noindent
{\bf Case I. }$K$ has $d+2$ vertices $x_1,\ldots,x_{d+2}$.

\vskip 1mm \noindent
{\bf (Ia) }There are $d+1$ vertices in one facet of $K$, say
$F=\conv(\{x_1,\ldots,x_{d+1}\})$. Then $K=\conv(\{x_{d+2}\}\cup F)$ is
a pyramid with basis $F$. Using Lemma~\ref{lem-G} and induction on the
dimension, we get the inequality, we also realize that equality here is
possible only if $K$ is a simplex, which it is not.

\vskip 1mm \noindent
{\bf (Ib) }No $d+1$ vertices of $K$ are in the same hyperplane. Thus $K$
is a simplicial polytope. Then there is a hyperplane $L$, containing $x_1$
(or any other vertex of $K$), and a relative neighborhood of $x_1$,
$U\subset L$, such that, for any $y\in U$, the polytope
$\conv(\{y,x_2,\ldots,x_{d+2}\})$ has the same volume as $K$ and
$y,x_2,\ldots,x_{d+2}$ are its vertices. For any non-zero vector $v$
parallel to the
hyperplane $L$, substituting $y=x_1+tv$ above, provides a volume preserving
shadow system $K_t$ for $t\in [-a,b]$ for some $a,b>0$. From
Theorem~\ref{th-A} it follows that at least in one direction, say for
$t>0$, $|K^*_t|$ is non-increasing. Since $K$ is not a pyramid
with apex $x_1$, we can find a direction $v$ as above, such that for $t=b$ the
moving point $x_1+tv$ hits the hyperplane of a facet of $K$ that does not
contain $x_1$. Thus $K_b$ is either a simplex or of the type (Ia) above. As
we have seen, $\Pi_d(K_t)$ is non-increasing for $t>0$ thus the inequality is
proved in Case (Ib) as well.
Assume now that $K$ is minimal, then
$|K^*_t|$ is also non-decreasing in both directions. It follows (again with
the help of Theorem~\ref{th-A}) that
$|K^*_t|$ must be constant throughout $[-a,b]$. Choosing the direction $v$
as above, we see that the two options at which we have reached for $K_b$
are impossible, because, by Proposition~\ref{pr-A}, $K_t$ are affine
images of $K$ for all $t\in [-a,b]$.

We conclude that Case I contains no minimal polytope.

\vskip 1mm \noindent
{\bf Case II. }$K$ has $d+3$ vertices $x_1,\ldots,x_{d+3}$.

\vskip 1mm \noindent
{\bf (IIa) }There are $d+2$ vertices in one facet of $K$. Then $K$ is a
pyramid and, like in the case (Ia), we get the inequality by induction
and $K$ can not be minimal.

\vskip 1mm \noindent
{\bf (IIb) }There are $d+1$ vertices of $K$ in one hyperplane, but no $d+2$
vertices are in one hyperplane. We may assume that $x_3,\ldots,x_{d+3}$ are
such that $F=\conv(\{x_3,\ldots,x_{d+3}\})\subset \rd{d-1}$ (again,
$\rd{d}$ is represented as $\rd{d-1}\times \R$) and that
$x_1=(X_1,\xi_1)$, $x_2=(X_2,\xi_2)$ satisfy $\xi_1\xi_2\neq 0$. We
distinguish here between three possible situations:

\vskip 1mm \noindent
{\bf (IIb1) }$\xi_1\xi_2<0$, say $\xi_1<0<\xi_2$. That is, $x_1$ and
$x_2$ are on opposite
sides of $\rd{d-1}$. Then, by the convexity of $K$, the line segment
$]x_1,x_2[$ meets $F$. $K$
is then a double pyramid with basis $F$. Let $v=x_2-x_1$ and define
a shadow system $K_t$ (with $K_0=K$) by keeping all the vertices of
$F$ fixed and moving $x_i$ to $x_i+tv$ for $i=1,2$. That is, we shift
the segment $[x_1,x_2]$ at constant speed along its line
and $K_t$ is the convex hull of
$F$ with the shifted segment. If $t\in ]\tau_1,\tau_2[=
]\frac{-\xi_2}{\xi_2-\xi_1},\frac{-\xi_1}{\xi_2-\xi_1}[$ then $K_t$ is
still a double pyramid and $|K_t|=|K|$. $K_{\tau_1}$ and $K_{\tau_2}$
are pyramids with basis $F$ and apex $x_1$, respectively $x_2$. By
Theorem~\ref{th-A}, $|K^*_t|$ is non-increasing in at least one of the
directions, say for $t>0$. But $K_{\tau_2}$ is of type (Ia) or (IIa).
This proves the inequality and also shows that $K$ is not minimal.

\vskip 1mm \noindent
{\bf (IIb2) }$\xi_1\xi_2>0$ and $\xi_1\neq\xi_2$, say $0<\xi_1<\xi_2$.
Let $x_0\in \rd{d}$ be the intersection point of the line $L$ through
$x_1$ and $x_2$ with $\rd{d-1}$. Clearly $x_0\not\in F$. Let
\[V_2=|\conv(\{x_0\}\cup F)|\mbox{ and } V_1=V_2-|F|\]
($(d-1)$-dimensional volumes). Then
\begin{equation}
\label{eq-G}
|K|=\frac{1}{d}(\xi_2V_2-\xi_1V_1)\,.
\end{equation}
Let $v=x_2-x_1$ and define the shadow system $K_t$ (with $K_0=K$) by
keeping the vertices of $F$ fixed, moving $x_1$ to
\[x_1(t)=(X_1(t),\xi_1(t))=x_1+tv\,,\]
and $x_2$ to
\[x_2(t)=(X_2(t),\xi_2(t))=x_2+t\frac{V_1}{V_2}v\,.\]
That is, $x_1$ and $x_2$ move along their joint
line in the same direction but at different speeds. If $t\in
[\tau_1,\tau_2]=[\frac{-\xi_1}{\xi_2-\xi_1},
\frac{1}{1-\frac{V_1}{V_2}}]$
then $0\leq \xi_1(t)\leq \xi_2(t)$ and we have, like in (\ref{eq-G})
\[|K_t|=\frac{1}{d}(\xi_2(t)V_2-\xi_1(t)V_1)=
\frac{1}{d}\left(\left(\xi_2+t\frac{V_1}{V_2}(\xi_2-\xi_1)\right)V_2
-(\xi_1+t(\xi_2-\xi_1))V_1\right)=|K|.\]

The bodies $K_{\tau_1}$ and $K_{\tau_2}$ are pyramids, $K_{\tau_1}$
has basis $\conv(\{x_0\}\cup F)$ and apex $x_2(\tau_1)$ and $K_{\tau_2}$
has basis $F$ and apex $x_2(\tau_2)=x_1(\tau_2)$. Since, by
Theorem~\ref{th-A}, $|K^*_t|$ is non-increasing at least in one direction,
we get the inequality by induction. Also $K$ can not be minimal because,
had it been minimal, $|K^*_t|$ would be constant in the interval
$[\tau_1,\tau_2]$ (we have seen the argument before), but $K_{\tau_2}$
is of type (Ia) and thus it is not minimal.

\vskip 1mm \noindent
{\bf (IIb3) }$\xi_1\xi_2>0$ and $\xi_1=\xi_2=\xi$, say $\xi>0$. That is,
the line $L$ connecting $x_1$ and $x_2$ is parallel to $\rd{d-1}$. The
inequality can be proved in this case by considering $K$ as the limit
(in the Hausdorff metric) of polytopes of type (IIb2). Since we wish
to characterize the minimal case, we give another proof. Let $P_L$ be the
orthogonal projection in $\rd{d-1}$ onto the $(d-2)$-dimensional subspace
orthogonal to $x_2-x_1$. Then
\begin{equation}
\label{eq-H}
|K|=\frac{1}{d}\xi\left(|F|+\frac{|x_2-x_1|}{n-1}|P_L(F)|\right)\,.
\end{equation}

We define a shadow system $K_t$ (with $K_0=K$) by fixing all the vertices
of $K$ except $x_2$ which moves (along the line $L$) to $x_2+tv$, where
$v=x_2-x_1$. For $t\in [-1,\infty[$ the formula analogous to (\ref{eq-H})
is
\begin{equation}
\label{eq-I}
|K_t|=|K|+t\xi\frac{|x_2-x_1|}{d(d-1)}|P_L(F)|\,.
\end{equation}
Note that $K_{-1}$ is a pyramid of the type (Ia) and that, as $t$ tends to
$\infty$, bounded affine images of $K_t$ converge in the Hausdorff metric
to a pyramid with basis which is the orthogonal projection of $K$ onto the
hyperplane orthogonal to $v$ (this pyramid may happen to be a simplex).

\begin{lem}
\label{lem-H}
Let $U$ be a open convex subset of $\rd{d}$, $\phi:U\rightarrow \R$
a positive convex function and $\psi:U\rightarrow \R$ a non-negative
concave function. Then $\frac{\psi}{\phi}$ does not attain its minimum
in $U$, unless $\frac{\psi}{\phi}$ is constant in $U$, in which case
$\psi=c\phi$ is affine.
\end{lem}

\vskip 1mm \noindent
{\bf Proof of Lemma~\ref{lem-H}. }If $a\in U$ satisfies
\[\frac{\psi(a)}{\phi(a)}\leq \frac{\psi(x)}{\phi(x)}\]
for every $x\in U$ then, setting
\[h(x)=\frac{\psi(a)}{\phi(a)}\phi(x)-\psi(x)\,,\]
we see that $h$ is a convex function that attains a maximum at $a$ and
$h(a)=0$. This is possible only if $h$ is identically zero. In this case
$\psi=\frac{\psi(a)}{\phi(a)}\phi$ is concave and convex, and thus affine.
$\frac{\psi}{\phi}=\frac{\psi(a)}{\phi(a)}$ is constant in $U$ in this case.
\bbx
\medskip

Applying Lemma~\ref{lem-H} to $\psi(t)=|K_t|$ (which, by (\ref{eq-I}), is
affine - hence concave) and
$\phi(t)=|K^*_t|^{-1}$, we conclude the following behavior of $\Pi_d(K_t)$:
if it is not constant in $[-1,\infty[$ then it either attains its minimum
at $t=-1$ or tends to its infimum as $t$ tends to $\infty$. $K_{-1}$, being
of type (Ia), is not minimal. Hence $K=K_0$ can not be minimal. Also, since
the limiting body, as $t$ tends to $\infty$, of bounded images of $K_t$,
is a pyramid, we conclude that $K$ satisfies the inequality of
Theorem~\ref{th-B}.

\vskip 1mm \noindent
{\bf (IIc) }No $d+1$ vertices of $K$ are in the same hyperplane. In this
case $K$ is simplicial. This case is treated in the same way as case
(Ib), to show the inequality and the fact that $K$ is not minimal.

We have thus checked all the possible configurations of the vertices of $K$
and verified that no polytope with more than $d+1$ vertices can be minimal.
\bbx

\section{New insight into known results}

We demonstrate in this section how the tools that were developed in Section~2,
and in particular the investigation of the case of equality, can be applied
to provide ``natural'' proofs of two known results.

The first one of these is the Blaschke-Santal\'o inequality \cite{Sant}. A
simple proof of a more general inequality is given in \cite{M-P}. The
proof of \cite{M-P} uses the inequality part of our Lemma~\ref{lem-A}, the
proof of this inequality as presented in Section 2 here, is taken
from there. Using
Proposition~\ref{pr-A}, we are able to give here a smooth form of
the proof.
Particularly simple here is the characterization of the maximal bodies as
ellipsoids. A step that in \cite{M-P} required reduction to the centrally
symmetric case and reference to a lemma of Saint-Raymond \cite{SR} (another
proof of the characterization of maximal bodies is in \cite{Pet}, that proof
requires deep results in PDE, together with a complicated reduction to the
smooth case).

\begin{theo}[Blaschke, Santal\'o, Saint-Raymond, Petty]
\label{th-C}
Let $K$ be a convex body in $\rd{d}$; then
\begin{equation}
\label{eq-J}
\Pi_d(K)\leq \Pi_d(B_2^d)\,,
\end{equation}
where $B_2^d$ is the Euclidean unit ball in $\rd{d}$. Equality holds in
(\ref{eq-J}) if and only if $K$ is an ellipsoid.
\end{theo}

\vskip 1mm \noindent
{\bf Proof. }We accept the fact, which is proved in a standard manner, that
$\max \Pi_d(K)$ is attained among the convex bodies $K\subset \rd{d}$. Let
$H$ be a hyperplane in $\rd{d}$ and $K_H$ the result of Steiner
symmetrization of $K$ about $H$. $K_H$ can be considered as $K_0$ of a shadow
system $K_t$, $t\in [-c,c]$, where $K_{-c}=K$ and $K_c$ is the mirror
reflection of $K$ about $H$. This shadow system preserves the lengths of
chords of $K_t$ that are orthogonal to $H$. Thus $|K_t|=|K|$ for all $t$.

Also, as $K_t$ is an affine image (reflection) of $K_{-t}$, it follows that
$|K_t^*|=|K_{-t}^*|$ for $t\in [-c,c]$. By Theorem~\ref{th-A} we conclude
that $|K_t^*|$ attains its maximum at $t=0$. Moreover, if $K$ is a convex body
at which $\max \Pi_d(K)$ is attained, it follows from Proposition~\ref{pr-A}
that $K_H$ is an affine image of $K$. That is, the midpoints of all the
chords of $K$ that are orthogonal to $H$ lie in a hyperplane. This happen
for {\em any\/} hyperplane $H$. We conclude, using a classical (basically
2-dimensional) result of Brunn
(see \cite{Blasch}), that $K$ is an ellipsoid. \bbx
\vskip 2mm

The second result that we treat is reverse Santal\'o inequality in dimension 2
(without symmetry assumption). Mahler \cite{Mah} proved that among polygons in
$\rd{2}$, the minimum of $\Pi_2(K)$ is achieved only by triangles. Clearly the
proof for polygons proves the inequality for general convex bodies in $\rd{2}$.
The case of equality, however, does not follow. The characterization of
triangles as the only convex bodies in $\rd{2}$ that are minimal for $\Pi_2$
has been given by Meyer \cite{Mey-2}. The proof in  \cite{Mey-2} is tricky
(see \cite{Mey-Rei} for a generalization of the method and some applications).
Here we present a ``natural'' proof of the inequality together with the case
of equality.

The next lemma is known as a ``classical folklore''.

\begin{lem}
\label{lem-I}
Let $\cal C$ be the cone of concave, continuous functions on an interval
$[\alpha,\beta]$, that satisfy $f(\alpha)=f(\beta)=0$. The extreme (non-zero)
rays of $\cal C$ are spanned by the functions $f\in {\cal C}$ such that, for
some $\gamma\in ]\alpha,\beta[$, $f$ is affine on $[\alpha,\gamma]$ and on
$[\gamma,\beta]$. Denoting the set of these functions by $\cal R$, the above
claim means that if $f\in {\cal C}\setminus {\cal R}$ then there exist
$g,h\in {\cal C}$, both not proportional to $f$, such that $f=g+h$.
\end{lem}

\vskip 1mm \noindent
{\bf Proof. }We may assume that the interval $[\alpha,\beta]$ is $[0,1]$.
If $f\in {\cal C}$ and $a\in ]0,1[$ define
\[g(x)=f(x)-x(f(a)+(1-a)f'_L(a)) \mbox{\ for $x\in [0,a]$ and\ }\]
\[g(x)=(1-x)(f(a)-af'_L(a)) \mbox{\ for $x\in ]a,1]$.}\]
where $f'_L(a)$ is the left derivative of $f$ at $a$ ($f'_R$ is the right
derivative). Let $h=f-g$.

It is easy to verify that $g$ is continuous at $a$ and $g'_L(a)=g'_R(a)$.
As $g=f-(\mbox{\em an affine function})$ on $[0,a]$, $g$ is affine on
$[a,1]$, and $g(0)=g(1)=0$, it follows that $g\in {\cal C}$. Similarly,
$h$ is in $\cal C$ (one checks easily that $h'_L(a)\geq h'_R(a)$). Clearly,
if $f\not\in {\cal R}$ than for any $a\in ]0,1[$ neither $g$ nor $h$ is
proportional to $f$. And yet $f=g+h$.

Let us remark that, beside the above direct proof, one can prove the lemma
by observing that any $f\in {\cal C}$ can be represented by a unique
positive measure $\mu$ on $[0,1]$, as
\[f(x)=\int_0^1G(x,y)\,d\mu(y)\,,\]
where $G(x,y)=\min((1-y)x,y(1-x))$. Thus the extreme rays are associated
with the Dirac measures. That is, $f\in {\cal C}$ spans an extreme ray if
and only if $f(x)=\lambda \min((1-a)x,a(1-x))$ for some $\lambda>0$ and
$a\in ]0,1[$.\bbx
\medskip

\begin{theo}[Mahler, Meyer]
\label{th-D}
Let $K$ be a convex body in $\R^{2}$. Then $\Pi_2(K)$
is minimal if and only if $K$ is a triangle.
\end{theo}

\vskip 2mm\noindent
{\bf Proof.}
Assume, without loss of generality, that a diameter of $K$ coincides with
the $x$-axis. Then
$$K=\{(x,y); x\in [\alpha,\beta], y\in [a(x),b(x)]\}\ ,$$
where, with the notations of the previous lemma,
the functions $-a$ and $b$ are in the cone $\cal C$ (the fact that
$a(\alpha)=a(\beta)=b(\alpha)=b(\beta)=0$ is due to the fact that the
$x$-axis contains a diameter of $K$).

Without loss of generality, we can assume that $b$ is not identically $0$.
Suppose
that $b= {b_{0}+b_{1}\over 2}$, with $b_{0}, b_{1}\in {\cal C}$ and
 both $b_{0}$ and $b_{1}$ different from $b$.

We define a shadow movement based on $M= \{(x,a(x)), (x,b(x)); x\in
[\alpha,\beta]\}$ with direction the $y$-axis, $t\in [-1,1]$
and speed defined by
$v(x,a(x))=0$ and $v(x, b(x)) ={{b_{1}(x)-b_{0}(x)\over 2}}$.
We then have $K_0= K$,\ $K_{-1}=\{(x,y); x\in [\alpha,\beta], y\in
[a(x),b_{0}(x)]\}$
 and $
 K_{1}=\{(x,y); x\in [\alpha,\beta], y\in [a(x),b_{1}(x)]\}. $
 Clearly
 \[t\mapsto |K_{t}|=|K|+ \frac{t}{2}\int_{\alpha}^{\beta}\big(
 b_{1}(x)-b_{0}(x) \big) dx\]
  is affine on $[-1,1]$.

If $K$ has minimal volume product (we again accept the standard
fact of the existence of a minimal body), it follows from
Proposition~\ref{pr-A}  that there exists an affine map
$T:\R^2\mapsto \R^2$
such that $T(K)=K_{1}$ and $T$ is of the form
$T(x,y)=(x,ux+vy+w)$ for some $u,v,w\in \R$.

By the Remark made after the proof of Proposition~\ref{pr-A} we have
for every $x\in[\alpha,\beta]$: $a(x)=va(x) + ux+w$ and
$b_1(x)=vb(x) + ux+w$.

Since $a$ vanishes at $\alpha$ and $\beta$, one has $u=w=0$. Now,
since $b_{1} \not=b$, one has $v\not=1$ and thus $a\equiv 0$ and $b_{1}=vb$
(observe that $v\neq 0$ because $b_{1}\neq 0$).

We have thus shown that $b$ generates an extreme ray of ${\cal C}$ and that
$a\equiv 0$. Hence, by Lemma~\ref{lem-I}, $K$ is a triangle.  \bbx

\end{document}